\documentclass[11pt]{article}
\usepackage{amsmath, amsthm, amssymb}    
\usepackage{newtxtext}
 
\usepackage[margin=2.5cm,noheadfoot]{geometry}
\usepackage{graphicx}			
\usepackage[colorlinks=true, urlcolor=blue, citecolor=black, linkcolor=black]{hyperref}  
\usepackage{gensymb}
\usepackage{float}
\usepackage{graphicx,wrapfig,lipsum}
\usepackage{subcaption}			
\usepackage{xfrac}
\usepackage[symbol]{footmisc}
\usepackage{sectsty}
\sectionfont{\centering}

\renewenvironment{thebibliography}[1]{
  \begin{oldthebibliography}{#1}
    \setlength{\itemsep}{0em}
    \setlength{\parskip}{0em}
}
{
  \end{oldthebibliography}
}

\title{Exploring Mathematics with Curvagon Tiles}

\author{Hanne Kekkonen\textsuperscript{} 
\vspace{10pt}\\
\textsuperscript{}EEMCS, Delft University of Technology, Netherlands; h.n.kekkonen@tudelft.nl} 

\date{}					

\begin{document}

\maketitle

\thispagestyle{empty}

\begin{abstract}
\footnotesize
\noindent Building blocks and tiles are an excellent way of learning about geometry and mathematics in general. 
There are several versions of tiles that are either snapped together or connected with magnets that can be used to introduce topics like volume, tessellations, and Platonic solids. However, since these tiles are made of hard plastic, they are not very suitable for creating hyperbolic surfaces or shapes where the tiles need to bend. 
Curvagons are flexible regular polygon building blocks that allow you to quickly build anything from hyperbolic surfaces and tori to dinosaurs and shoes. They can be used to introduce mathematical concepts from Archimedean solids to Gauss-Bonnet theorem. You can also let your imagination run free and build whatever comes to mind. 
\end{abstract}

\section{Introduction}\label{Sec:Intro}

\begin{figure}[H]
\centering
\begin{minipage}[b]{0.9\textwidth} 
	\includegraphics[width=\textwidth]{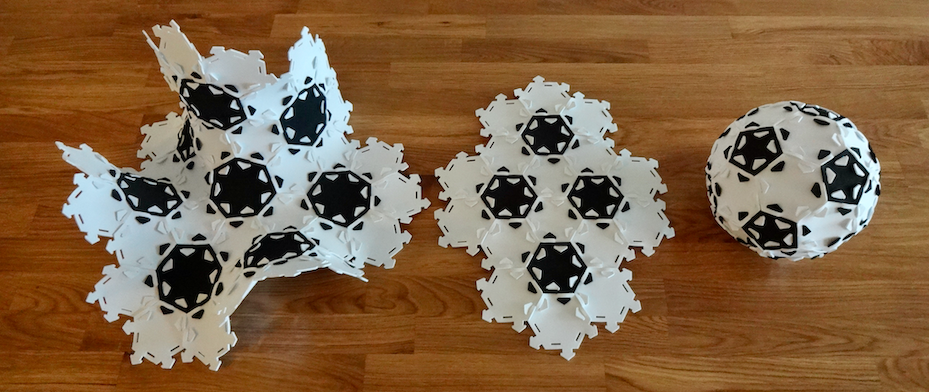}
\end{minipage}
\caption{`Football' tiling models of the hyperbolic plane, Euclidean plane, and sphere using Curvagons.}\label{Fig:HyperbolicFootball}
\label{Fig:Tiling}
\end{figure}

Introducing new and surprising concepts to students is a good way to get them excited about mathematics. One of these concepts is non-Euclidean geometry. Since we live on a globe, spherical geometry is quite intuitive to many but it also offers a couple of surprises because most people have only encountered Euclidean geometry in school. An even more interesting topic is hyperbolic geometry of surfaces with a constant negative curvature. Both of these concepts can be introduced using polygon tilings. One way to approximate a sphere is to start with a pentagon and surround it with five hexagons. When continued, this pattern creates the traditional football approximation of the sphere (truncated icosahedron), as shown in Figure \ref{Fig:HyperbolicFootball}. To create an approximation of a hyperbolic plane, where every point is a saddle point, we can start with a heptagon instead of a pentagon and surround it with seven hexagons.
This `hyperbolic football' model has the advantage that, while the model is curved, each polygon is flat and so you can draw straight lines on it. For more details see \cite{Henderson, Sottile}. 
In the model introduced by Keith Henderson, the polygons are printed on paper, cut, and then taped together.
I really liked the idea and tested the model in a couple of hands-on outreach events, but noticed that many students struggled to build a sufficiently large model during the time that was reserved for the activity. 

\noindent Inspired by the paper version of the hyperbolic football, I wanted to create building blocks that could be used to quickly create surfaces with negative curvature. The building blocks should be flexible allowing the surface to bend naturally and the connection mechanism should be symmetric since half of the polygons have odd number of sides. The result is Curvagons; flexible regular polygon tiles (all angles and edges identical) that are interlocked together to create a supple surface that can be bent and twisted smoothly.
These tiles can be made with a cutting machine from different materials depending on the intended use. Tiles made from EVA (ethylene-vinyl acetate) craft foam can be used in a classroom to build everything from hyperbolic surfaces to dinosaurs. You can also use narrow tape to create straight lines, which allows further exploration of hyperbolic geometry, as shown in Figure \ref{Fig:curvatures}.  
Tear-resistant paper and paper-like materials, e.g. leather paper, can be used to create durable models that can be drawn on with a normal pencil, as shown in Figure \ref{Fig:Triangles}.


\section{Mathematics with Curvagons}
\subsection{Non-Euclidean Geometry}
\begin{wrapfigure}{r}{5.1cm}
\includegraphics[width=5.1cm]{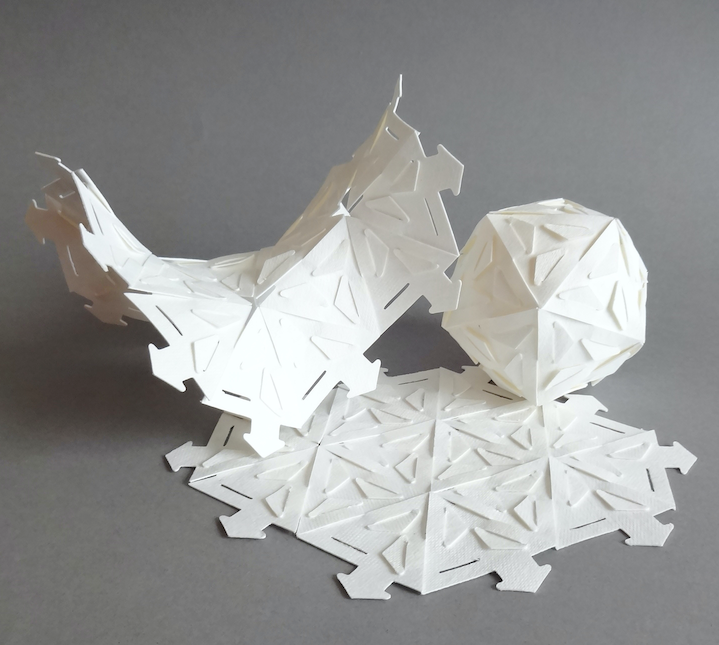}
\caption{Tilings made from SnapPap Curvagons.}
\label{Fig:Triangles}
\end{wrapfigure}

Consider six equilateral triangles meeting at a vertex. Since the angles of a triangle add up to $180\degree$ in Euclidean geometry the internal angles of an equilateral triangle are $60\degree$ and six of them add up to $360\degree$ (or equivalently $2\pi$ radians) creating a flat surface, as shown in Figure \ref{Fig:Triangles}. If you instead place five triangles at a vertex, the surface will curve towards itself, and if continued, close creating a crude approximation of a sphere called icosahedron. The amount by which the sum of the angles at a vertex falls short of $360\degree$ is called angle defect. For the icosahedron, this shortfall is $60\degree$ since one triangle is missing. We can also approximate the hyperbolic plane by putting seven equilateral triangles at a vertex, creating a $\{3,7\}$ polyhedral model. This model of the hyperbolic plane is rather wrinkly due to the $-60\degree$ angle defect which we from now on call $60\degree$ angle excess. The hyperbolic plane looks locally like the Euclidean plane meaning that we can get smoother approximations of it using tessellations where the angle excess is small. One of these models is the `hyperbolic football' tessellation described in the Introduction, where the angle excess in every vertex is only $8\sfrac{4}{7}\degree$. 

\begin{figure}[H]
\centering
\begin{minipage}[b]{0.9\textwidth} 
	\includegraphics[width=\textwidth]{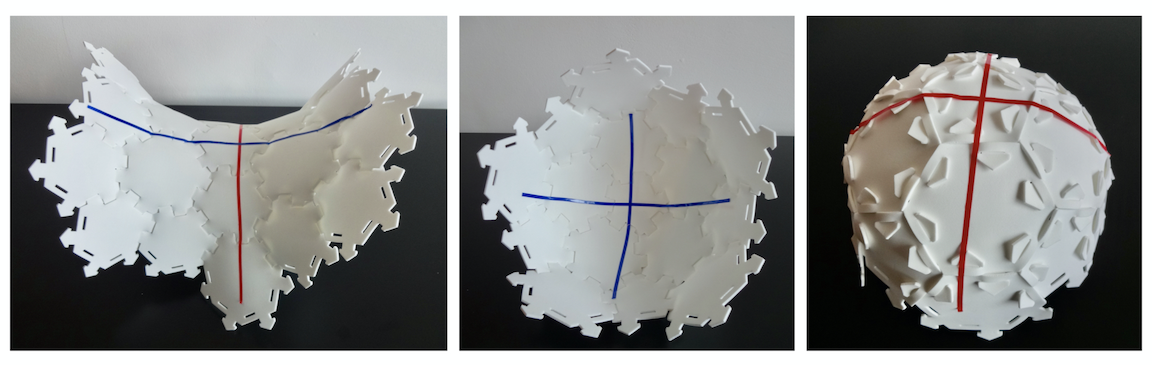}
\end{minipage}
\caption{One of the largest curves on a saddle shape looks like a hill (red) while the other one looks like a valley (blue). For a sphere both of the curves are either hills or valleys.}
\label{Fig:curvatures}
\end{figure}

Notice that while the plane is flat, the sphere and the hyperbolic plane curve. Looking at the model of the hyperbolic plane, we see that it creates a saddle shape and we can find orthogonal directions with largest curving. One of these curves looks like a hill while the other one looks like a valley, as shown in Figure \ref{Fig:curvatures}. To distinguish hills from valleys we can consider one of them as having positive curvature and the other\linebreak
\vspace{-4mm}

\noindent  one as having negative curvature. We call these maximum and minimum values of curvature the principal curvatures. If we do the same for the sphere, we notice that both of the lines curve in the same direction and are the same size, therefore the maximum and minimum curvature are the same. 
The Gaussian curvature, that describes how a surface curves at a point, is given as the product of the two principal curvatures. So the hyperbolic plane has negative Gaussian curvature while spheres have positive Gaussian curvature. 
The curvature of a circle is simply given by $1/R$, where $R$ is the radius of the circle, and so the curvature of a sphere is $1/R^2$.
This means that larger spheres have smaller curvature, as we would expect. Similarly, the milder the saddle shape of the hyperbolic plane, the closer the curvature is to zero. 

\begin{wrapfigure}{l}{6cm}
\includegraphics[width=6cm]{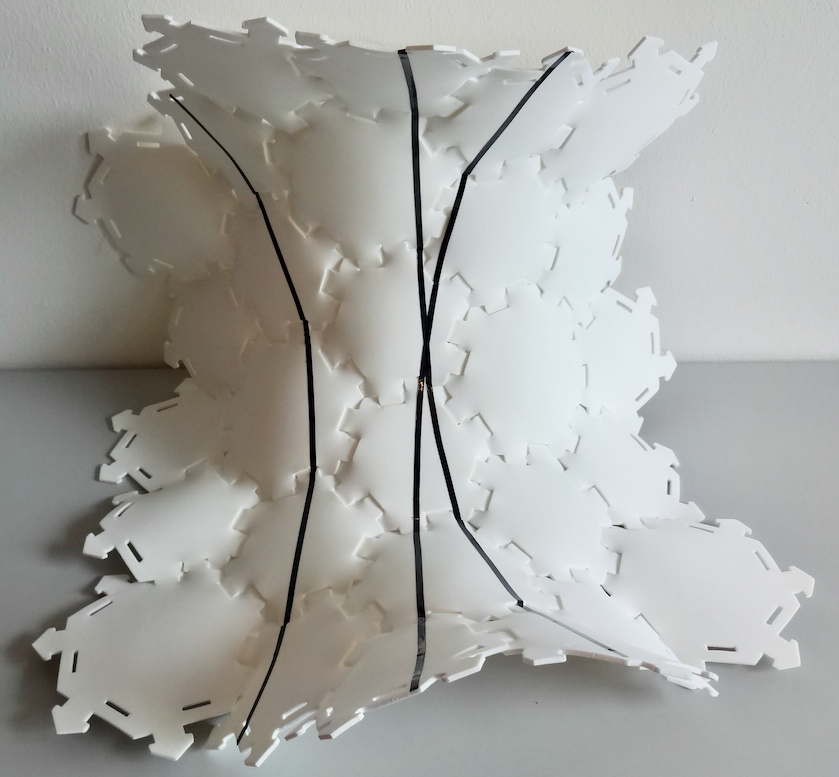}
\caption{Diverging parallels on an EVA foam Curvagon model.}
\label{Fig:Lines}
\end{wrapfigure}  
Many of the fundamental truths you learn in school are not true when you leave the Euclidean plane. For example, on a hyperbolic surface the angles of a triangle add up to less than $180\degree$, a line has infinitely many parallel lines through a given point, and the circumference of a circle is more than $2\pi R$. Also, even though  a hyperbolic surface is larger than the flat plane you cannot draw arbitrarily large triangles on them. This might seem unintuitive at first but becomes clear if you study how parallel lines behave on hyperbolic surfaces. You can draw a straight line on an EVA foam `hyperbolic football' model by gently straightening it along a line with a ruler and tracing the line using a narrow tape.
You can then straighten the model along another line that does not cross the original line. You will notice that the lines are closest to each other at one point but diverge in both directions, as shown in Figure \ref{Fig:Lines}. If you start drawing a very large triangle you will notice that this is not possible because the lines never meet. They might seem to converge at first but after reaching a certain point, where they are closest to each other, they start diverging again. The largest triangle you can draw on a hyperbolic plane is called the ideal triangle and the sum of its angles is zero.   
Notice that even though the ideal triangle has finite area its endpoints stretch to infinity so drawing one would be rather challenging. For other models introducing hyperbolic geometry see \cite{Kekkonen} and for more in-depth introduction to curvature and hyperbolic geometry see e.g. \cite{Taimina1, Taimina2}. You can find other hyperbolic geometry activities for the hyperbolic football from \cite{Sottile}.

\subsection{Gauss-Bonnet Theorem}
\begin{wrapfigure}{r}{5cm}
\includegraphics[width=5cm]{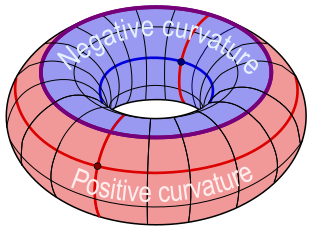}
\caption{A torus has regions with positive (red), negative (blue) and zero curvature (purple).}
\label{Fig:TorusCurvature}
\end{wrapfigure}  
Spheres and hyperbolic surfaces have the same Gaussian curvature everywhere, which is why we say that they have constant curvature. Since they look the same at every point they can be easily approximated using simple tessellations as described above. But how about less regular shapes? We start by considering a torus. If you put your finger on the outer edge of a torus, you will notice that the surface curves away from your hand, much like a sphere does, as shown in Figure \ref{Fig:TorusCurvature}. These points on the outer edge have positive Gaussian curvature. If you put your finger on the inner circle of the torus, the surface curves away from your hand in one direction but in the other direction it curves towards it, creating a saddle shape. This region has negative Gaussian curvature. If you place the torus on a table, you can see that on the very top and bottom, between the regions of positive and negative curvature, there is a circle that lies perfectly flat and hence has zero curvature. So, to approximate a torus we need to create a tessellation with angle excess in the middle and angle defect on the outer edge. But how much excess and defect should we have?

\begin{wrapfigure}{r}{5.7cm}
\includegraphics[width=5.7cm]{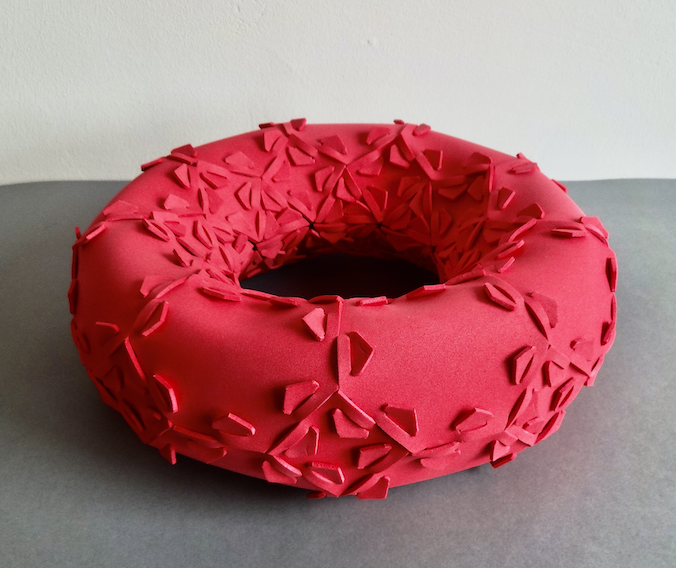}
\caption{An approximation of a torus}
\label{Fig:TorusTesselation}
\end{wrapfigure}  
Here we can consider several more complex mathematical concepts with advanced students. One of these is the Euler characteristic which describes the shape or structure of a topological space regardless of the way it is bent. We first consider shapes without a boundary like the torus and sphere. The Euler characteristic is then given by $\chi=vertices-edges+faces$ which was later proven to equal to $\chi=2-2*genus$. Genus is the number of doughnut-like holes on the surface. Thus a sphere has genus zero, a torus (and a coffee cup) has genus one, and a double torus has genus two. Since a torus has genus one, its Euler characteristic is zero.
The approximation of a torus in Figure \ref{Fig:TorusTesselation} consists of $2*9$ triangles, $3*9$ squares and $2*9$ heptagons, and it has $9*9$ vertices, $16*9$ edges and $7*9$ faces, which means that its Euler characteristic given by the first definition is indeed zero.  

We can also compute the angle defects and excess' of all the vertices of the torus in Figure \ref{Fig:TorusTesselation} and notice that they add up to zero. This means that the angle defect on the outer surface is as large as the angle excess in the middle. The result is an example of the Descartes' theorem on the total defect of a polyhedron, which states that the sum of the total defect is $2\pi*\chi$. 
 Descartes' theorem is a special case of the Gauss-Bonnet theorem, where the curvature is concentrated at discrete points, the vertices.
The Gauss-Bonnet theorem is a fundamental result that bridges the gap between differential geometry and topology by connecting the geometrical concept of curvature with the topological concept of Euler characteristic. It states that for closed surfaces the total curvature equals to $2\pi*\chi$. 
We can conclude that since the Euler characteristic of a torus is zero, its total curvature is also zero. 
In general the Gauss-Bonnet theorem means that no matter how you bend a surface (without creating holes) its total curvature will stay the same since its Euler characteristic, being a topological invariant, does not change even though the Gaussian curvature at some points does. 
For example, if you push a dimple into a sphere, the total curvature of the surface will stay $4\pi$ even though the Gaussian curvature changes around the dimple.
Notice that Curvagons are made of solid flat pieces so the curvature is indeed concentrated at the vertices. Other way to approximate curved shapes is to create pieces with holes, like Curvahedra \cite{Harriss}, where the curvature is spread over the open region, or specially design pieces where the curvature is pushed into meandering edges as in Zippergons \cite{Delp}.


We can also consider hyperbolic triangles to see what extra information the Gauss-Bonnet theorem provides.  Calculating the Euler characteristic is easy since all triangles have three vertices, three edges, and one face, that is, $\chi=3-3+1=1$. Notice that triangles have a boundary and so they are quite different from the torus considered above. Because of the boundary, the Gauss-Bonnet theorem becomes a bit more 
 \linebreak
\vspace{-3mm}
 \begin{figure}[H]
\centering
\begin{minipage}[b]{0.71\textwidth} 
	\includegraphics[width=\textwidth]{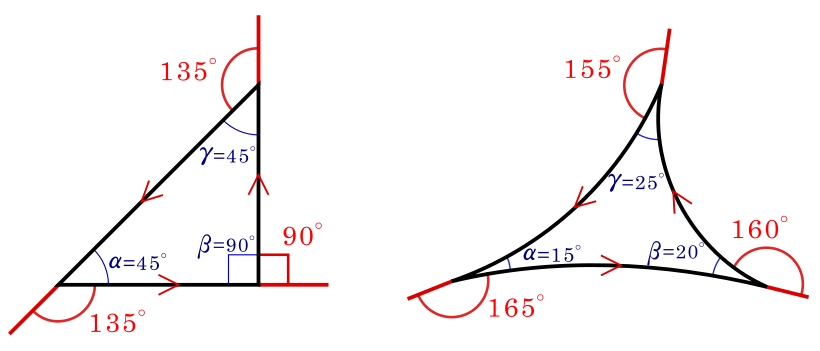}
\end{minipage}
\caption{Turning angles in red and the interior angles in blue. The total turning is the sum of the turning angles. For a flat triangle, this is $360\degree$. For a hyperbolic triangle, this is strictly more than $360\degree$. }
\label{Fig:TotalTurning}
\end{figure}
\begin{wrapfigure}{l}{5.23cm}
\includegraphics[width=5.23cm]{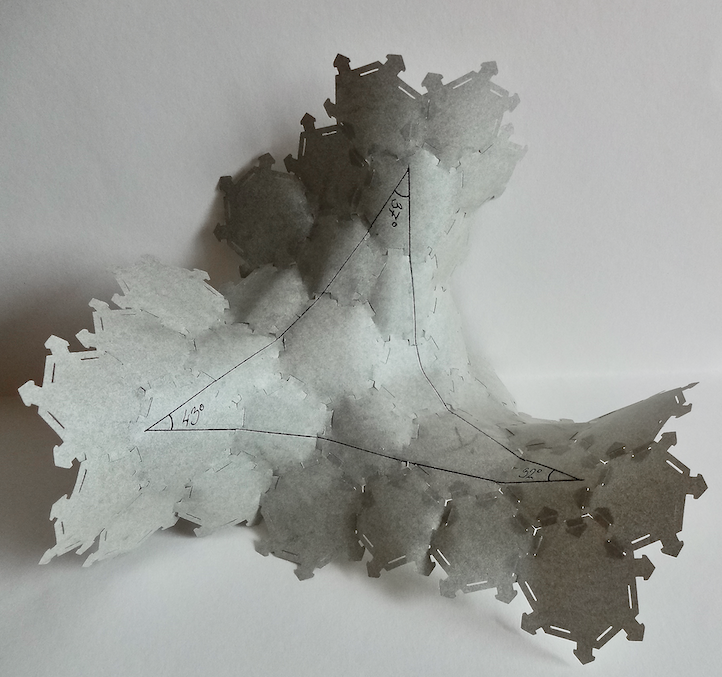}
\caption{A hyperbolic triangle on a Curvagon paper model.}
\label{Fig:HyperbolicTriangles}
\end{wrapfigure} 
\noindent 
complex stating that $total\ curvature\ within\ a \ triangle =2\pi*\chi- total\ turning =2\pi- total\ turning$. Total turning describes how much a person following the edge of a triangle would have to turn in order to 
return to the point they started from, as shown in Figure \ref{Fig:TotalTurning}. It is given by the sum of the turning angles, where
$turning\ angle=\pi-interior\ angle$.  Hence the total turning for a triangle can be written as $3\pi-\angle\alpha-\angle\beta-\angle\gamma$, where $\angle\alpha,\angle\beta,\angle\gamma$ are the interior angles of the triangle, and so we can write $total\ curvature =2\pi- (3\pi-\angle\alpha-\angle\beta-\angle\gamma)=\angle\alpha+\angle\beta+\angle\gamma-\pi$.
We can now use the knowledge that the Euclidean plane has zero curvature to deduce the fact that the angles of a triangle add up to $\pi$ on a  flat surface. For hyperbolic triangles the total curvature within a triangle equals to $\angle\alpha+\angle\beta+\angle\gamma-\pi$, the deviation of the angle sum from the Euclidean triangle angle sum $\pi$. Since the curvature of the tiling model is concentrated at the vertices, the deviation should equal $-\ number\ of \ vertices\ inside\ the\ triangle*angle\ excess$. You can test this by drawing several triangles of different sizes on your model and calculating the angle sums. 
Figure \ref{Fig:HyperbolicTriangles} shows a triangle on a hyperbolic football model. It has has eight vertices inside it, so the curvature within the triangle is $-8*8\sfrac{4}{7}\degree\approx -68.57\degree$. The interior angles measure $43\degree$, $32\degree$ and $37\degree$, and hence the deviation of the angle sum  from $180\degree$ is $-68\degree$. Notice that the largest triangle you could in theory draw on an infinitely large  hyperbolic football model would have $21$ vertices inside it and all the angles would be zero. 

 \subsection{Platonic and Archimedean Solids}

\begin{figure}[H]
\centering
\begin{minipage}[b]{0.98\textwidth} 
	\includegraphics[width=\textwidth]{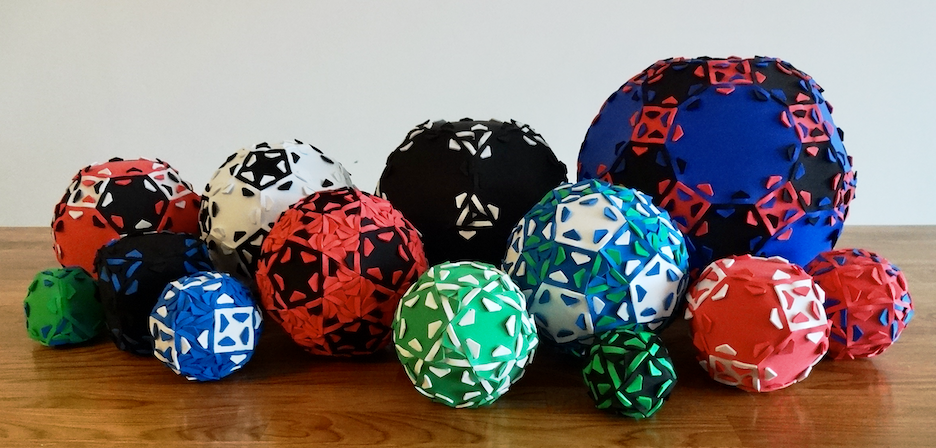}
\end{minipage}
\caption{The $13$ Archimedean solids build from EVA foam Curvagons.}
\label{Fig:Solids}
\end{figure}

Platonic solids are convex regular polyhedra, which means that all their faces are identical regular polygons and that the same number of faces meet at each vertex. You can easily test that there exists only five platonic solids. Notice that a vertex needs at least 3 faces and an angle defect. If the angle defect is zero the regular tiling will fill the Euclidean plane. If there is angle excess you will get a saddle shape. 
Next we relax the above rules and allow the use of different regular polygons to create convex polyhedra called Archimedean solids, as shown in Figure \ref{Fig:Solids}. All the vertices are still expected 
\linebreak
\vspace{-1mm}

\begin{wrapfigure}{l}{6cm}
\includegraphics[width=6cm]{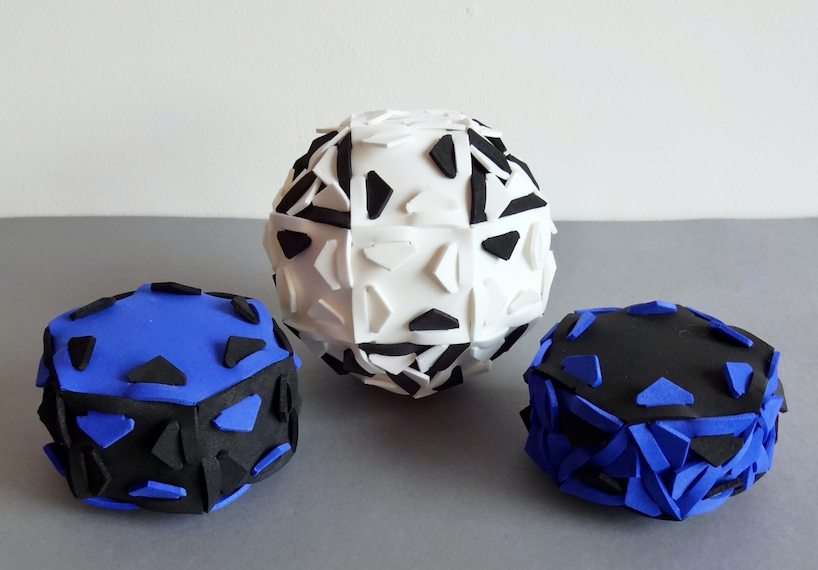}
\caption{A prism (left), elongated square gyrobicupola (middle), and antiprism (right) are only locally symmetric.}
\label{Fig:Prisms}
\end{wrapfigure}   
\noindent  
to be identical and the shapes have to be globally symmetric. The symmetry requirement means that 
prisms, antiprisms and the elongated square gyrobicupola are not considered to be Archimedean solids,
as shown in Figure \ref{Fig:Prisms}. 
There are 13 Archimedean solids (excluding the five Platonic solids) and they can be constructed from the Platonic solids by truncation, that is, by cutting away corners. 

Notice that all of these convex polyhedra can be considered to be approximations of the sphere. They have no holes in them so their genus is zero and so their Euler characteristic $\chi$ is two. Using the Descartes' theorem we see that the total angle defect of all of the solids is $2\pi*\chi=4\pi$ or $720\degree$. This can be quite surprising since the smallest solid (tetrahedron) has only four vertices while the biggest solid (truncated icosidodecahedron) has $120$. The result  is explained by the fact that the more vertices a solid has the smaller the angle defects become. The angle defects of a tetrahedron are $180\degree$ while the angle defects of a truncated icosidodecahedron are only $6\degree$.
Since all the vertices of Platonic and Archimedean solids are identical, the total defect is given by $angle\ defect*number\ of\ vertices$ which means that $number\ of\ vertices=720\degree/angle\ defect$.


\section{Polygon Sculpting}

Approximating shapes and understanding curvature is important in many applications. For example, polygon, and especially triangle, meshes are used in 3D computer graphics. Also, to create a well-fitting dress you need to know how to cut patterns to create positive and negative curvature. Since Curvagons are flexible and can be woven together quickly they can be used for testing new ideas and for polygon sculpting. As a simple exercise, students can consider what kind of packaging and capsules can be created using different polygons, and which of them can be opened in a way that would result in near zero-waste cutting patterns. 

We can also consider how to approximate more complex shapes using Curvagons. Regular polygon tessellations have long been used to create crocheted blankets but these so-called granny squares (or, more generally, polygons) have also been used to make stuffed animals. We can use Curvagons instead of crocheted
\begin{figure}[H]
\centering
\begin{minipage}[b]{0.86\textwidth} 
	\includegraphics[width=\textwidth]{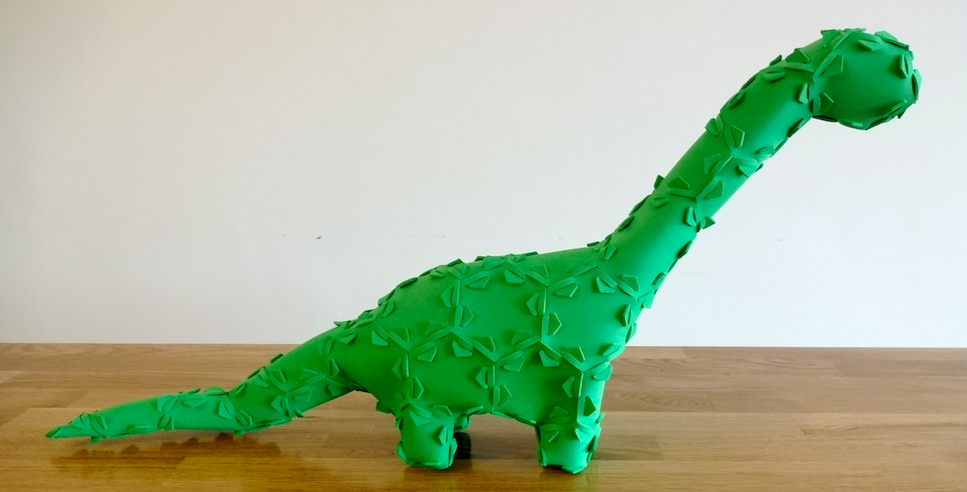}
\end{minipage}
\caption{A polygon brontosaurus made from EVA foam Curvagons.}
\label{Fig:Bronto}
\end{figure}

\begin{wrapfigure}{r}{6.7cm}
\includegraphics[width=6.7cm]{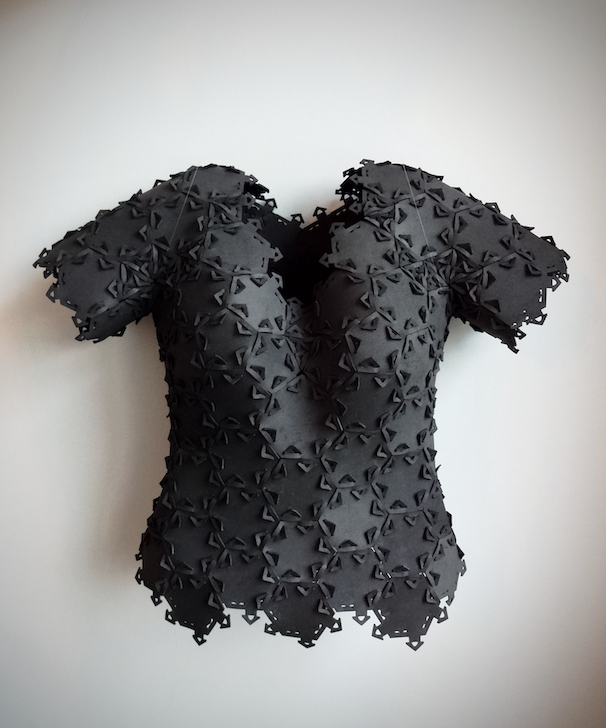}
\caption{A polygon approximation of a torso from EVA foam Curvagons.}
\label{Fig:Torso}
\end{wrapfigure}  
\noindent polygons to test different animal models. See Figure \ref{Fig:Bronto} for a Curvagon brontosaurus. 
Building different animals allows students to get creative while still forcing them to think about how different areas curve and how to mimic this using polygons. Polygonal animals can also be used as an introduction to computer graphics where smooth shapes are estimated by triangle meshes. How would you need to divide the regular polygons of your model to create a good piecewise flat  approximation?

\begin{figure}[H]
\centering
\begin{minipage}[b]{0.93\textwidth} 
	\includegraphics[width=\textwidth]{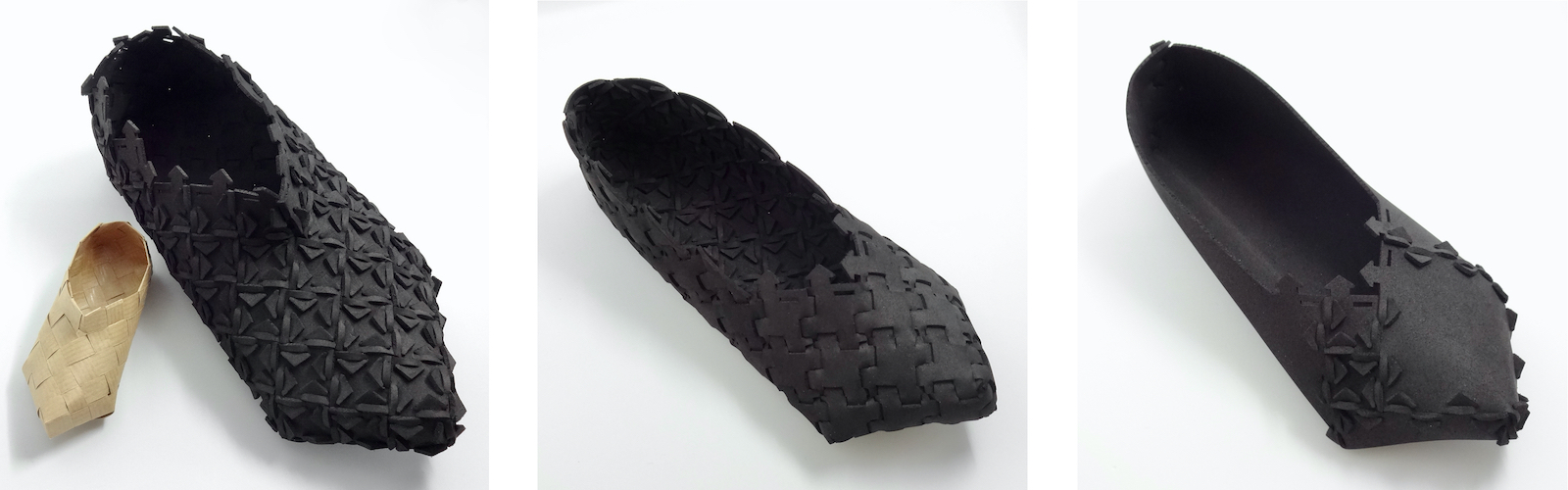}
\end{minipage}
\caption{A traditionally woven {virsu} and the same design using Curvagons (left), a more slipper-like design showing the interesting pattern on the smooth side (middle), and a {virsu} inspired slipper (right).}
\label{Fig:Shoes}
\end{figure}

Weaving has long been used to create clothing and other utility articles. One example is the Finnish \textit{virsu}, a type of bast shoe traditionally woven from strips of birch bark. The weave creates a square tiling where corners can be made by placing three squares around a vertex and a saddle shape (for the ankle) by using five squares. One can use small square Curvagons to create different versions of \textit{virsu}, as shown in Figure \ref{Fig:Shoes}. These models can also be opened several ways to find cutting patterns for \textit{virsu} type slippers. Curvagon slippers could be made from leather or felt. Since the individual pieces are small, one can use discarded scrap materials which would otherwise go to waste. 
These type of modular techniques, where small pieces are used to create larger surfaces, have gained popularity in recent years by repurposing waste materials from fashion houses and tanneries to create sustainable fashion. 
Interlocking small square Curvagons creates an interesting texture on the smooth side of the model. This simple weave could be used to create pillowcases, rugs, bags and many other everyday items from different leftover materials.

We can also take the idea of modular fashion a bit further and 
use Curvagons to study the different types of shapes needed for creating clothes in general. A well-fitting shirt has areas of positive, negative and zero curvature, as shown in Figure \ref{Fig:Torso}. 
A tight pencil skirt takes much less fabric but requires more modelling than a 50s-style skirt, which can be laid out to create a full circle, and so could be constructed using an Euclidean tiling. The inspiration between mathematics and fashion can also work to the other direction as was illustrated by Zippergons \cite{Delp}. Zippergons are flexible patterned pieces that are designed to fit a given shape and were inspired by discussions between a mathematician Bill Thurston and Dai Fujiwara, the director of design for the Issey Miyake fashion house, on how designing fitting clothes and mathematics is related. 


 \newpage
 \section{Summary and Conclusions}
 
\begin{figure}[H]
\centering
\begin{minipage}[b]{0.84\textwidth} 
	\includegraphics[width=\textwidth]{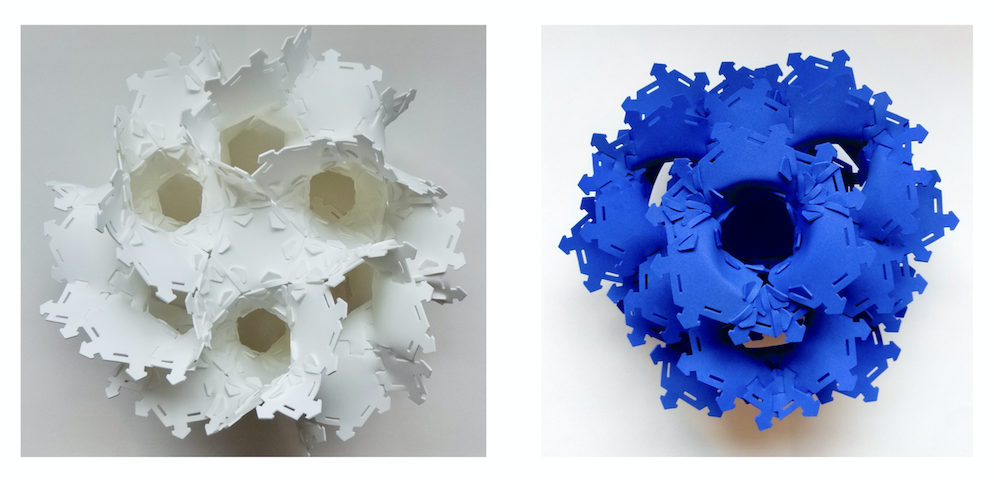}
\end{minipage}
\caption{A Curvagon approximation of the Schoen's Gyroid (left) and the Schwarz' D surface (right), which are examples of triply periodic minimal surfaces.}
\label{Fig:MinimalSurfaces}
\end{figure}

This paper demonstrates how Curvagon tiles can be used to introduce and explore several mathematical concepts. Curvagons  are flexible regular polygon building blocks that can be made from different materials and they can be assembled into infinitely many possible shapes. They can be used to introduce mathematical concepts suitable for different educational levels ranging from elementary school to university. Younger students can experiment with angles, different ways to tile the plane, and the Platonic and Archimedean solids. Older students can be introduced to the concept of non-Euclidean geometry and how many geometry facts they have learned in school are actually true only in Euclidean geometry. They can also test how all the Platonic and Archimedean solids (or any polyhedron that is homeomorphic to the sphere) have the same total angle defect. Hyperbolic geometry is also interesting for university-level students and Curvagons can even be used to introduce fundamental results like Gauss-Bonnet theorem through its special case of Descartes' theorem on total angle defect.  Curvagons can also be used to construct beautiful mathematical objects like triply periodic minimal surfaces such as Schoen's Gyroid and the Schwarz' D surface, as shown in Figure \ref{Fig:MinimalSurfaces}. In conclusion, Curvagons can be used to construct just about any shape you can imagine.

{\setlength{\baselineskip}{12pt} 
\raggedright				

} 

\end{document}